%% file: main.tex
\documentclass[11pt]{article}

\usepackage{sectsty}
\usepackage{graphicx}
\usepackage{authblk}
\usepackage{hyperref}
\usepackage{xcolor}
\usepackage{mathtools}
\usepackage{amsfonts}
\usepackage{subcaption} % subcaptions for subfigures
\usepackage{thmtools}
\usepackage{thm-restate}

\usepackage{amsthm}
\usepackage{amssymb}

\newtheorem{theorem}{Theorem}

\newtheorem{assumption}{Assumptions}
\newtheorem{lemma}{Lemma}

\title{A Proof of Strong Consistency of Maximum Likelihood Estimator for Independent Non-Identically Distributed Data}
\author[1]{Ricardo Ferreira}
\author[1]{Filipa Valdeira}
\author[2]{Marta Guimar\~{a}es}
\author[1]{Cl\'{a}udia Soares}
\affil[1]{NOVA School of Science and Technology}
\affil[2]{Neuraspace}

\date{}

\begin{document}
\maketitle

\let\thefootnote\relax\footnotetext{Ricardo Ferreira, Filipa Valdeira and Cl\'{a}udia Soares are with NOVA School of Science and Technology
(e-mail: \href{rjn.ferreira@campus.fct.unl.pt}{rjn.ferreira@campus.fct.unl.pt},~\href{f.valdeira@fct.unl.pt}{f.valdeira@fct.unl.pt},~\href{claudia.soares@fct.unl.pt}{claudia.soares@fct.unl.pt}). Marta Guimar\~{a}es is with Neuraspace (e-mail: \href{marta.guimaraes@neuraspace.com}{marta.guimaraes@neuraspace.com}).}

% Optional TOC
% \tableofcontents
% \pagebreak

%--Paper--

\begin{abstract}
    We give a general proof of the strong consistency of the Maximum Likelihood Estimator for the case of independent non-identically distributed (i.n.i.d) data, assuming that the density functions of the random variables follow a particular set of assumptions. Our proof is based on the works of Wald~\cite{wald1949note}, Goel~\cite{goel1974note}, and Ferguson~\cite{ferguson2017course}. We use this result to prove the strong consistency of a Maximum Likelihood Estimator for Orbit Determination.
\end{abstract}

\input{intro}
\input{strong-consistency}

\input{mle-consistency-proof}

%--/Paper--

\bibliographystyle{ieeetr}
\bibliography{bibliography}

\end{document}

%% file: intro.tex
\section{Introduction}

Maximum Likelihood Estimation is an estimation method that has been intensively studied since presented by Fisher~\cite{hald1999history}. The goal is to estimate a set of parameters such that the joint probability distribution of the observed data is maximized.

Consider a set of $n$ random variables $X_i$, for $i = 1, \dots, n$, with a joint distribution $g(x_1, \dots, x_n;\theta)$ governed by some parameter vector $\theta \in \Theta$. So, the objective of Maximum Likelihood Estimation is to find the parameter $\theta$ that maximizes the likelihood function

\begin{align*}
    \mathcal{L}(\theta) = g(x_1, \dots, x_n;\theta) .
\end{align*}
When the random variables are consider to be independent, then the likelihood function becomes the product of the individual density functions $g_i(x_i;\theta)$, i.e.,

\begin{align*}
    \mathcal{L}(\theta) = \prod_{i=1}^{n} g_i(x_i;\theta) .
\end{align*}
In practical terms, due to the fact that the logarithm function is a monotonic function, the optimization usually is performed considering the log-likelihood function

\begin{align*}
    \log \mathcal{L}(\theta) = \sum_{i=1}^{n} \log g_i(x_i;\theta),
\end{align*}
as it simplifies the computation in many practical cases~\cite{bishop2006pattern}.

Over the yeares, different works have proven different properties for the Maximum Likelihood Estimator under different circumstances, such as asymptotical normality and consistency of the estimator~\cite{ferguson2017course,rao1973linear}. In this work, we present a proof for the strong consistency of the Maximum Likelihood Estimator for the case of independent non-identically distributed random variables, based on the works of Wald~\cite{wald1949note}, Goel~\cite{goel1974note}, and Ferguson~\cite{ferguson2017course}. With this result, we show that the Maximum Likelihood Estimator used for orbit determination in~\cite{ferreira2024generalizing} is strongly consistent.

%% file: strong-consistency.tex
\section{Strong consistency of Maximum Likelihood Estimator for i.n.i.d random variables}

Firstly, in a generic way, for $i=1, 2, \dots, n$, consider the observations $x_i$ which are realizations of the i.n.i.d random variables $X_i$ with density $g_i(x;\theta)$ with respect to a $\sigma$-finite measure $\nu$, respectively. Assume that the parameter $\theta \in \Theta \subseteq \mathbb{R}^{d}$ and its true value $\theta_0 \in \Theta$. All the expected values, variances and probability statements are determined under the assumption that $\theta_0$ is the true parameter point. All the assumptions, lemmas, theorems and proofs are based on the works of Wald~\cite{wald1949note}, Goel~\cite{goel1974note}, and Ferguson~\cite{ferguson2017course}.

Now, assume that all density functions $g_i$ satisfy the set of assumptions presented in Assumptions~\ref{as:strong-consistency-assumptions}.

\begin{assumption}
    \label{as:strong-consistency-assumptions}
    \hfill
    \begin{itemize}
        \item[(i)] $\Theta$ is a compact set; \\
        \item[(ii)] (Identifiability) $g_i(x;\theta) = g_i(x;\theta_0)$ almost everywhere (a.e.) $\implies \theta = \theta_0$; \\
        \item[(iii)] The function $g_i(x;\theta, \rho) = \underset{\theta' \in \Theta}{\sup}\{g_i(x;\theta') : \|\theta' - \theta\| < \rho \}$, for $i = 1, 2, \dots$, is a measurable function of $x$ for any fixed $\theta \in \Theta$ and $\rho > 0$; \\
        \item[(iv)] For each $\theta \neq \theta_0$, each sufficiently small value of $\rho > 0$, and for some non-decreasing sequence $b_n$ of constants there exists an $\varepsilon > 0$ and an integer $N^{*}(\theta, \rho)$ such that for all $n \geq N^{*}$
        \begin{align*}
            \frac{1}{b_n} \sum_{i=1}^{n} \mathbf{E}\left[ \log g_i(X_i;\theta, \rho) - \log g_i(X_i;\theta_0) \right] < - \varepsilon ;
        \end{align*}

        \item[(v)] For each $\theta \neq \theta_0$, and any sufficiently small values of $\rho$,
        \begin{align*}
            \sum_{i=1}^{\infty} \frac{1}{b_i^2} \mathbf{Var}\left[ \log g_i(X_i;\theta, \rho) - \log g_i(X_i;\theta_0) \right] < \infty .
        \end{align*}

    \end{itemize} 
\end{assumption}
Based on these assumptions, we consider the next lemma, which will help to prove our first theorem.

\begin{lemma}
    \label{lemma:lemma-1}
    For any $\theta \neq \theta_0$, we have
    \begin{align*}
        \mathbf{E}\left[ \log g_i(X_i;\theta) \right] < \mathbf{E}\left[ \log g_i(X_i;\theta_0) \right].
    \end{align*}
\end{lemma}

\begin{proof}
    Due to Assumption (ii), this lemma follows from the Shannon-Kolmogorov information inequality~\cite{ferguson2017course}.
\end{proof}

\begin{theorem}
\label{th:theorem-1}
    If $X_1, \dots, X_n$ are i.n.i.d random variables whose density functions satisfy Assumptions~\ref{as:strong-consistency-assumptions}, and $\omega$ is any closed subset of the parameter set $\Theta$ which does not contain the true parameter $\theta_0$, then
    \begin{align*}
        \mathbf{P} \left( \underset{n \to \infty}{\lim} \frac{\underset{\theta \in \omega}{\sup} \prod_{i=1}^{n} g_i(X_i; \theta)}{\prod_{i=1}^{n} g_i(X_i; \theta_0)} = 0 \right) = 1.
    \end{align*}
\end{theorem}

\begin{proof}
    Define $\varphi_{i}(x,r) \coloneqq \underset{\|\theta\| > r}{\sup} g_i(x;\theta)$. Fix $\varepsilon > 0$ and let $r_0$ be such that
    \begin{equation}
        \centering
        \label{eq:eq-1}
        \begin{aligned}
            \mathbf{E}\left[ \log \varphi_{i}(X; r_0) \right] < \mathbf{E}\left[ \log g_i(X;\theta_0) \right] - \varepsilon , \quad \mbox{for } i = 1, 2, \dots, n,
        \end{aligned}
    \end{equation}
    and let $\omega_1 = \{ \theta : \theta \in \omega, \|\theta\| \leq r_0 \}$. For each point $\theta \in \omega_1$, assume there is a $\rho_{\theta} > 0$, such that $\mathbf{E}\left[ \log g_i(X;\theta, \rho_{\theta}) \right] < \mathbf{E}\left[ \log g_i(X;\theta_0) \right]$, for $i = 1, 2, \dots, n$. The existence of $r_0$ and $\rho_\theta$ derive from Lemma~\ref{lemma:lemma-1} and Assumption (i) that $\Theta$ is compact. Since $\omega_1$ is a closed subset of $\Theta$, it is also a compact set, thus there exists a finite number of points $\theta_1, \dots, \theta_h \in \omega_1$ such that $B(\theta_1, \rho_{\theta_1}) + \dots + B(\theta_h, \rho_{\theta_h})$ contains $\omega_1$ as a subset, where $B(\theta, \rho_{\theta})$ denotes the ball of center $\theta$ and radius $\rho_{\theta}$. So, we have
        \begin{align*}
            0 &\leq \underset{\theta \in \omega}{\sup} \prod_{i=1}^{n} g_i(X_i; \theta) \\
            &\leq \sum_{j=1}^{h} \left[ \prod_{i=1}^{n} g_i(X;\theta_j, \rho_{\theta_j}) + \prod_{i=1}^{n} \varphi_{i}(X;r_0) \right] .
        \end{align*}
    Therefore, we can prove Theorem~\ref{th:theorem-1} if we show that 
    \begin{align*}
        \mathbf{P} \left( \underset{n \to \infty}{\lim} \frac{\prod_{i=1}^{n} g_i(X_i; \theta_j, \rho_{\theta_j})}{\prod_{i=1}^{n} g_i(X_i; \theta_0)} = 0 \right) = 1, \quad \mbox{for } j = 1, \dots, h, \quad \mbox{and} \quad \mathbf{P} \left( \underset{n \to \infty}{\lim} \frac{\prod_{i=1}^{n} \varphi_i(X_i; r_0)}{\prod_{i=1}^{n} g_i(X_i; \theta_0)} = 0 \right) = 1 .
    \end{align*}
    We can rewrite these equations as
    \begin{equation}
        \centering
        \label{eq:eq-2}
        \begin{aligned}
            \mathbf{P} \left( \underset{n \to \infty}{\lim} \sum_{i=1}^{n} \left[ \log g_i(X_i; \theta_j, \rho_{\theta_j}) - \log g_i(X_i; \theta_0) \right] = -\infty \right) = 1, \quad \mbox{for } j = 1, \dots, h,
        \end{aligned}
    \end{equation}
    and
    \begin{equation}
        \centering
        \label{eq:eq-3}
        \begin{aligned}
            \mathbf{P} \left( \underset{n \to \infty}{\lim} \sum_{i=1}^{n} \left[ \log \varphi_{i}(X; r_0) - \log g_i(X_i; \theta_0) \right] = -\infty \right) = 1.
        \end{aligned}
    \end{equation}
    Both Eqs.~\eqref{eq:eq-2} and~\eqref{eq:eq-3} can be easily obtained by combining Assumptions (iv) and (v), with Eq.~\eqref{eq:eq-1} and Kolmogorov's strong law of large numbers~\cite{sen1994large}. Thus, we prove Theorem~\ref{th:theorem-1}.
\end{proof}
With this result, we can now show that the MLE estimator for i.n.i.d random variables whose density functions follow the set of Assumptions~\ref{as:strong-consistency-assumptions} is strongly consistent.

\begin{theorem}
    \label{th:theorem-2}
    Let $\Bar{\theta}_{n}(x_1, \dots, x_n)$ be a function of the observations $x_1, \dots, x_n$, which are realizations of the i.n.i.d random variables $X_i$, whose density functions satisfy Assumptions~\ref{as:strong-consistency-assumptions}, such that
    \begin{align*}
        \frac{\prod_{i=1}^{n} g_i(x_i; \Bar{\theta}_n)}{\prod_{i=1}^{n} g_i(x_i; \theta_0)} \geq c > 0,
    \end{align*}
    for all $n$ and all $(x_1, \dots, x_n)$. Then $\mathbf{P}\left(\underset{n \to \infty}{\lim} \Bar{\theta}_n = \theta_0 \right) = 1$.
\end{theorem}

\begin{proof}
    We must show that for any $\varepsilon > 0$, the probability is one that all limit points $\Bar{\theta}$ of the sequence $\{\Bar{\theta}_n\}$ satisfy the inequality $\|\Bar{\theta} - \theta_0\| \leq \varepsilon$. Assuming that there exists a limit point $\Bar{\theta}$ such that $\|\Bar{\theta} - \theta_0\| > \varepsilon$ implies that 
    \begin{align*}
        \underset{\|\theta - \theta_0\| > \varepsilon}{\sup} \prod_{i=1}^{n} g_i(x_i;\theta) \geq \prod_{i=1}^{n} g_i(x_i;\Bar{\theta}_n),
    \end{align*}
    for infinitely many $n$. However, then
    \begin{align*}
        \frac{\underset{\|\theta - \theta_0\| > \varepsilon}{\sup} \prod_{i=1}^{n} g_i(x_i; \theta)}{\prod_{i=1}^{n} g_i(x_i; \theta_0)} \geq c > 0,
    \end{align*}
    for infinitely many $n$. Since, according to Theorem~\ref{th:theorem-1}, this is an event with probability zero, then we proved that with probability one, all the limit points $\Bar{\theta}$ of the sequence $\{\Bar{\theta}_n\}$ satisfy the inequality $\|\Bar{\theta} - \theta_0\| \leq \varepsilon$.
\end{proof}
The MLE estimator, if it exists, satisfies Theorem~\ref{th:theorem-2} with $c = 1$. Therefore we establish the strong consistency of the MLE estimator for i.n.i.d random variables whose density functions satisfy Assumptions~\ref{as:strong-consistency-assumptions}. So, to show the strong consistency of our MLE estimator, we must show that these assumptions are satisfied.

%% file: mle-consistency-proof.tex
\section{Strong consistency of Maximum Likelihood Estimator for Orbit Determination}

As an example, we will apply this result to the estimator in~\cite{ferreira2024generalizing}, by showing that satisfies Assumptions~\ref{as:strong-consistency-assumptions}. As a summary of the work in~\cite{ferreira2024generalizing}, the goal is to determine a six-parameter state vector that describes the motion of a satellite, specifically the position, $\mathbf{r} \in {\mathbb{R}}^{3}$, and velocity, $\mathbf{v} \in {\mathbb{R}}^{3}$, in the Cartesian reference frame.

We consider a set of $N$ monostatic radars (capable of transmitting and receiving its own signal), located at $\mathbf{s}_i \in {\mathbb{R}}^{3}$, for $i = 1, \ldots, N$, all able to observe the satellite at the same instant. We assume that each radar $i$ is able to collect a set of different measurements:

\begin{itemize}
    \item Range measurements, $d_i = \|\mathbf{r} - \mathbf{s}_i\| + \epsilon_{d_i}$;
    \item Angle measurements, $ \mathbf{u}_i = \frac{\mathbf{r} - \mathbf{s}_i}{\|\mathbf{r} - \mathbf{s}_i\|} + \epsilon_{\mathbf{u}_i}$;
    \item Doppler shift measurements
    \begin{align*}
        f_i = \frac{2 f_{c,i}}{c} \left(\frac{\mathbf{r} - \mathbf{s}_i}{\|\mathbf{r} - \mathbf{s}_i\|}\right)^T \mathbf{v} + \epsilon_{f_i}.
    \end{align*}
\end{itemize}
where $c$ denotes the speed of light and $f_{c,i}$ denotes the carrier frequency of the signal from radar $i$. Range and Doppler shift noise is modeled as a Gaussian distribution with zero mean and variance $\sigma^2_{d_i}$ and $\sigma^2_{f_i}$, respectively, so $\epsilon_{d_i} \sim \mathcal{N}\left(0, \sigma^2_{d_i}\right)$ and $\epsilon_{f_i} \sim \mathcal{N}\left(0, \sigma^2_{f_i}\right)$.

Angle measurements are represented as unit-norm vectors and can be retrieved from elevation and azimuth measurements. Therefore, we model angle noise as a von Mises-Fisher distribution with mean direction zero and concentration parameter $\kappa_i$, which can better represent the distribution of directional data, so $\epsilon_{\mathbf{u}_i} \sim \mathcal{VMF}\left( \mathbf{0}, \kappa_i \right)$, $\mathbf{0} \coloneqq (0,0,0)$.

By assuming that the noise is independent, we can obtain the maximum likelihood estimator by solving the optimization problem

\begin{equation}
    \label{eq:mle-estimator-problem}
    \centering
    \begin{aligned}
    \underset{\mathbf{r}, \mathbf{v}}{\mbox{minimize}} \quad & f_{\text{range}}(\mathbf{r}) + f_{\text{angle}}(\mathbf{r}) + f_{\text{doppler}}(\mathbf{r},\mathbf{v})
    \end{aligned}
\end{equation}
such that the functions are defined as

\begin{equation}
    \label{eq:f-range-function}
    \centering
    \begin{aligned}
    f_{\text{range}}(\mathbf{r}) = \sum_{i=1}^{N} \frac{1}{2 \sigma_{d_i}^2} \left( \|\mathbf{r} - \mathbf{s}_i\| - d_i \right)^2 ,
    \end{aligned}
\end{equation}
\begin{equation}
    \label{eq:f-angle-function}
    \centering
    \begin{aligned}
    f_{\text{angle}}(\mathbf{r}) = -\sum_{i=1}^{N} \kappa_i \mathbf{u}_i^T \frac{\mathbf{r} - \mathbf{s}_i}{\|\mathbf{r} - \mathbf{s}_i\|} ,
    \end{aligned}
\end{equation}
\begin{equation}
    \label{eq:f-doppler-function}
    \centering
    \begin{aligned}
    f_{\text{doppler}}(\mathbf{r},\mathbf{v}) = \sum_{i=1}^{N} \frac{1}{2 \sigma_{f_i}^2} \left( \frac{2 f_{c,i}}{c} \left(\frac{\mathbf{r} - \mathbf{s}_i}{\|\mathbf{r} - \mathbf{s}_i\|}\right)^T \mathbf{v} - f_i \right)^2,
    \end{aligned}
\end{equation}
for the range, angle and Doppler shift measurements, respectively.

\subsection{Fulfillment of the necessary assumptions}

Before we demonstrate that our estimator fulfills the necessary assumptions to be strongly consistent, we reiterate the practical assumption that each radar is able to obtain a tuple of three types of measurements, $x_i \coloneqq (d_i, \mathbf{u_i}, f_i)$, for $i = 1, \dots, N$. This practical assumption will be important to prove the identifiability assumption.

\begin{itemize}
    \item[(i)] $\Theta$ is a compact set; \\

    The first assumption is easy to demonstrate as the parameters have some physical constraints. Remember that in our case, the parameter we are trying to estimate is the object orbital state in the Cartesian reference frame, i.e., $\theta \coloneqq (\mathbf{r}, \mathbf{v})$. In the extreme, the speed of the object must not surpass the speed of light, thus $\| \mathbf{v} \| \leq c$, where $c$ in this case denotes the speed of light. On the other hand, the position parameter $\mathbf{r}$ is the position of a satellite on an Earth-centered referential so we can bound it to a volume as big as we want, i.e., we can choose a constant $M$ such that $\| \mathbf{r} \| \leq M$. Also, we must consider the object to be above the Earth's surface, so we can also consider a constant $m$ greater than the radius of the Earth, such that $\| \mathbf{r} \| \geq m$. From these constraints, we see that $\Theta$ is a closed and bounded set, therefore, by the Heine-Borel theorem~\cite{rudin1964principles}, $\Theta$ is a compact set. \\
    
    \item[(ii)] (Identifiability) $g_i(x;\theta) = g_i(x;\theta_0)$ almost everywhere (a.e.) $\implies \theta = \theta_0$; \\

    Assume that the random variables $X_i$ describe a tuple of the three types of measurements obtained by each radar, i.e., $x_i \coloneqq (d_i, \mathbf{u_i}, f_i)$, for $i = 1, \dots, N$. Remember that we consider $N$ radars, each one able to obtain a measurement of each type. We assume that $d_i \sim \mathcal{N}\left( \|\mathbf{r} - \mathbf{s_i} \|, \sigma_{d_i}^2 \right)$, $\mathbf{u_i} \sim \mathcal{VMF}\left( \frac{\mathbf{r} - \mathbf{s}_i}{\|\mathbf{r} - \mathbf{s}_i\|}, \kappa_i \right)$ and $f_i \sim \mathcal{N}\left( \frac{2 f_{c,i}}{c} \left(\frac{\mathbf{r} - \mathbf{s}_i}{\|\mathbf{r} - \mathbf{s}_i\|}\right)^T \mathbf{v}, \sigma_{f_i}^2 \right)$. Since we assume each measurement is independent from each other, then we can see the density function of the tuple as $g_i(x_i;\theta) = g_{d_i}(d_i;\theta) \ g_{u_i}(\mathbf{u}_i;\theta) \ g_{f_i}(f_i;\theta)$, such that $g_{d_i}(x_i;\theta) = \mathcal{N}\left( d_i; \|\mathbf{r} - \mathbf{s_i} \|, \sigma_{d_i}^2 \right)$, $g_{u_i}(\mathbf{u}_i;\theta) = \mathcal{VMF}\left( u_i; \frac{\mathbf{r} - \mathbf{s}_i}{\|\mathbf{r} - \mathbf{s}_i\|}, \kappa_i \right)$ and $g_{f_i}(f_i;\theta) = \mathcal{N}\left(f_i ; \frac{2 f_{c,i}}{c} \left(\frac{\mathbf{r} - \mathbf{s}_i}{\|\mathbf{r} - \mathbf{s}_i\|}\right)^T \mathbf{v}, \sigma_{f_i}^2 \right)$. So, to prove identifiability we must show that 
    \begin{align*}
        g_{d_i}(d_i;\theta) = g_{d_i}(d_i;\theta_0) \bigwedge g_{u_i}(u_i;\theta) = g_{u_i}(u_i;\theta_0) \bigwedge g_{f_i}(f_i;\theta) = g_{f_i}(f_i;\theta_0) \implies \mathbf{r} = \mathbf{r}_0 \bigwedge \mathbf{v} = \mathbf{v}_0 .
    \end{align*} 
    Letting $\mathbf{z} = \mathbf{r} - \mathbf{s_i}$ and  $\mathbf{z}_0 = \mathbf{r}_0 - \mathbf{s_i}$, first we start with the density function of the range measurements.
    \begin{align*}
        g_{d_i}(d_i;\theta) = g_{d_i}(d_i;\theta_0) &\iff \frac{1}{\sqrt{2 \pi \sigma_{d_i}^2}} \exp\left( - \frac{(d_i - \|\mathbf{z}\| )^2}{2 \sigma_{d_i}^2} \right) = \frac{1}{\sqrt{2 \pi \sigma_{d_i}^2}} \exp\left( - \frac{(d_i - \|\mathbf{z}_0\| )^2}{2 \sigma_{d_i}^2} \right) \\
        &\iff (d_i - \|\mathbf{z}\|)^2 = (d_i - \|\mathbf{z}_0\| )^2 \\
        &\iff \|\mathbf{z}\|^2 - \|\mathbf{z}_0\|^2 - 2 d_i \|\mathbf{z}\| + 2 d_i \|\mathbf{z}_0\| = 0 \\
        &\iff (\|\mathbf{z}\| + \|\mathbf{z}_0\|) (\|\mathbf{z}\| - \|\mathbf{z}_0\|) - 2 d_i (\|\mathbf{z}\| - \|\mathbf{z}_0\|) = 0 \\
        &\iff \|\mathbf{z}\| - \|\mathbf{z}_0\| = 0 \ \vee \ \|\mathbf{z}\| + \|\mathbf{z}_0\| - 2 d_i = 0 \\
        &\iff \|\mathbf{z}\| = \|\mathbf{z}_0\| \ \vee \ \|\mathbf{z}\| + \|\mathbf{z}_0\| = 2 d_i \\
        &\iff \|\mathbf{r} - \mathbf{s_i}\| = \|\mathbf{r}_0 - \mathbf{s_i}\| \ \vee \ \|\mathbf{r} - \mathbf{s_i}\| + \|\mathbf{r}_0 - \mathbf{s_i}\| = 2 d_i .
    \end{align*} 
    Now, analyzing the density functions for the angle measurements, we obtain
    \begin{align*}
        g_{u_i}(\mathbf{u}_i;\theta) = g_{u_i}(\mathbf{u}_i;\theta_0) &\iff \frac{\kappa_i}{4 \pi \sinh{\kappa_i}} \exp\left( \kappa_i \mathbf{u}_i^T \frac{\mathbf{z}}{\|\mathbf{z}\|} \right) = \frac{\kappa_i}{4 \pi \sinh{\kappa_i}} \exp\left( \kappa_i \mathbf{u}_i^T \frac{\mathbf{z}_0}{\|\mathbf{z}_0\|} \right) \\
        &\iff \kappa_i \mathbf{u}_i^T \frac{\mathbf{z}}{\|\mathbf{z}\|} = \kappa_i \mathbf{u}_i^T \frac{\mathbf{z}_0}{\|\mathbf{z}_0\|} \\
        &\implies \frac{\mathbf{z}}{\|\mathbf{z}\|} = \frac{\mathbf{z}_0}{\|\mathbf{z}_0\|} \quad\quad\quad\quad\quad\quad \left(\forall \, \mathbf{k} , \mathbf{k}^T \mathbf{x} = \mathbf{k}^T \mathbf{y} \implies \mathbf{x} = \mathbf{y}\right) \\
        &\iff \frac{\mathbf{r} - \mathbf{s_i}}{\|\mathbf{r} - \mathbf{s_i}\|} = \frac{\mathbf{r}_0 - \mathbf{s_i}}{\|\mathbf{r}_0 - \mathbf{s_i}\|} .
    \end{align*} 
    where $\sinh(\cdot)$ denotes the hyperbolic sine function. By combining this result with the previous one for the density function of the range measurements, we have
    \begin{align*}
        \left( \|\mathbf{r} - \mathbf{s_i}\| = \|\mathbf{r}_0 - \mathbf{s_i}\| \ \bigvee \ \|\mathbf{r} - \mathbf{s_i}\| + \|\mathbf{r}_0 - \mathbf{s_i}\| = 2 d_i \right) \bigwedge \left( \frac{\mathbf{r} - \mathbf{s_i}}{\|\mathbf{r} - \mathbf{s_i}\|} = \frac{\mathbf{r}_0 - \mathbf{s_i}}{\|\mathbf{r}_0 - \mathbf{s_i}\|} \right) \implies \mathbf{r} = \mathbf{r}_0 .
    \end{align*}
    Finally, for the density function of Doppler shift measurements, let $M_i = \frac{2 f_{c,i}}{c}$, $w = \left(\frac{\mathbf{r} - \mathbf{s_i}}{\|\mathbf{r} - \mathbf{s_i}\|}\right)^T \mathbf{v}$ and $w_0 = \left(\frac{\mathbf{r}_0 - \mathbf{s_i}}{\|\mathbf{r}_0 - \mathbf{s_i}\|}\right)^T \mathbf{v}_0$, then we have
    \begin{align*}
        g_{f_i}(f_i;\theta) = g_{f_i}(f_i;\theta_0) &\iff \frac{1}{\sqrt{2 \pi \sigma_{f_i}^2}} \exp\left( - \frac{\left( f_i - M_i w \right)^2}{2 \sigma_{f_i}^2} \right) = \frac{1}{\sqrt{2 \pi \sigma_{f_i}^2}} \exp\left( - \frac{\left( f_i - M_i w_0 \right)^2}{2 \sigma_{f_i}^2} \right) \\
        &\iff \left( f_i - M_i w \right)^2 = \left( f_i - M_i w_0 \right)^2 \\
        &\iff M_i^2 w^2 - M_i^2 w_0^2 - 2 f_i M_i w + 2 f_i M_i w_0  = 0 \\
        &\iff (M_i w - M_i w_0) (M_i w + M_i w_0) - 2 f_i (M_i w - M_i w_0) = 0 \\
        &\iff w = w_0 \ \vee \ w + w_0 = \frac{2 f_i}{M_i} \\
        &\iff \left(\frac{\mathbf{r} - \mathbf{s_i}}{\|\mathbf{r} - \mathbf{s_i}\|}\right)^T \mathbf{v} = \left(\frac{\mathbf{r}_0 - \mathbf{s_i}}{\|\mathbf{r}_0 - \mathbf{s_i}\|}\right)^T \mathbf{v}_0 \ \bigvee \\ 
        &\quad\quad\quad\quad \left(\frac{\mathbf{r} - \mathbf{s_i}}{\|\mathbf{r} - \mathbf{s_i}\|}\right)^T \mathbf{v} + \left(\frac{\mathbf{r}_0 - \mathbf{s_i}}{\|\mathbf{r}_0 - \mathbf{s_i}\|}\right)^T \mathbf{v}_0 = \frac{2 f_i}{M_i} .
    \end{align*}
    The combination of this result with the previous one implies that $\mathbf{r} = \mathbf{r}_0$ and with the fact that $\forall \, \mathbf{k}, \mathbf{k}^T \mathbf{x} = \mathbf{k}^T \mathbf{y} \implies \mathbf{x} = \mathbf{y}$, implies that $\mathbf{v} = \mathbf{v}_0$. So, we have demonstrated the identifiability assumption.
    \\
    
    \item[(iii)] The function $g_i(x;\theta, \rho) = \underset{\theta' \in \Theta}{\sup}\{g_i(x;\theta') : \|\theta' - \theta\| < \rho \}$, for $i = 1, 2, \dots$, is a measurable function of $x$ for any fixed $\theta \in \Theta$ and $\rho > 0$; \\

    For any fixed $\theta \in \Theta$ and $\rho > 0$, the functions $g_i(x;\theta) = g_{d_i}(d_i;\theta) \ g_{u_i}(\mathbf{u}_i;\theta) \ g_{f_i}(f_i;\theta)$ are the product of continuous functions, therefore they are continuous. Thus $g_i(x;\theta, \rho)$ are also continuous since they are the pointwise supremum of continuous functions. So, if $g_i(x;\theta, \rho)$ are continuous, then it is also measurable. \\
    
    \item[(iv)] For each $\theta \neq \theta_0$, each sufficiently small value of $\rho > 0$, and for some non-decreasing sequence $b_n$ of constants there exists an $\varepsilon > 0$ and an integer $N^{*}(\theta, \rho)$ such that for all $n \geq N^{*}$
    \begin{align*}
        \frac{1}{b_n} \sum_{i=1}^{n} \mathbf{E}\left[ \log g_i(X_i;\theta, \rho) - \log g_i(X_i;\theta_0) \right] < - \varepsilon ;
    \end{align*}
    First, choose $\delta_{\mathbf{r}}, \delta_{\mathbf{v}} > 0$ such that $\|\mathbf{r}' - \mathbf{r}\| \leq \delta_{\mathbf{r}}$, $\|\mathbf{v}' - \mathbf{v}\| \leq \delta_{\mathbf{v}}$ and $\|\mathbf{\theta}' - \mathbf{\theta}\| \leq \rho$, where $\mathbf{\theta} \coloneqq (\mathbf{r}, \mathbf{v})$. For the density functions $g_{d_i}(d_i;\theta)$, $g_{u_i}(\mathbf{u}_i;\theta)$ and $g_{f_i}(f_i;\theta)$, define the functions $g_{d_i}(d_i;\theta,\rho)$, $g_{u_i}(\mathbf{u}_i;\theta,\rho)$ and $g_{f_i}(f_i;\theta,\rho)$ as in the previous assumption. Note that
    \begin{align*}
        \mathbf{E}\left[ \log g_i(X_i;\theta, \rho) - \log g_i(X_i;\theta_0) \right] &= \mathbf{E}\left[ \log g_{d_i}(D_i;\theta, \rho) - \log g_{d_i}(D_i;\theta_0) \right] \\
        &+ \mathbf{E}\left[ \log g_{u_i}(U_i;\theta, \rho) - \log g_{u_i}(U_i;\theta_0) \right] \\
        &+ \mathbf{E}\left[ \log g_{f_i}(F_i;\theta, \rho) - \log g_{f_i}(F_i;\theta_0) \right] ,
    \end{align*}
    where $D_i$, $U_i$, $F_i$ denote the random variables of the range, angle and Doppler shift measurements, respectively. As before, let $\mathbf{z} = \mathbf{r} - \mathbf{s}_i$ and $\mathbf{z}_0 = \mathbf{r}_0 - \mathbf{s}_i$. Given the constraints on $\delta$ and $\rho$, by the triangle inequality let $\mathbf{r} = \mathbf{r}_0 + \mathbf{\Delta_{r}}$, such that $\|\mathbf{\Delta_{r}}\| \leq \delta_{\mathbf{r}}$, and $\mathbf{v} = \mathbf{v}_0 + \mathbf{\Delta_{v}}$, such that $\|\mathbf{\Delta_{v}}\| \leq \delta_{\mathbf{v}}$. Let $\delta = \max\{\delta_{\mathbf{r}}, \delta_{\mathbf{v}}\}$. So, for the first term, we obtain
    \begin{align*}
        \mathbf{E}\left[ \log g_{d_i}(D_i;\theta, \rho) - \log g_{d_i}(D_i;\theta_0) \right] &= \mathbf{E}\left[ \frac{(D_i - \|\mathbf{z}_0\|)^2 - (D_i - \|\mathbf{z}\|)^2}{2 \sigma_{d_i}^2} \right] \\
        &= \mathbf{E}\left[ \frac{\|\mathbf{z}_0\|^2 - 2 D_i \|\mathbf{z}_0\| + 2 D_i \|\mathbf{z}\| - \|\mathbf{z}\|^2}{2 \sigma_{d_i}^2} \right] \\
        &= -\frac{(\|\mathbf{z}\| - \|\mathbf{z}_0\|)^2}{2 \sigma_{d_i}^2} \leq - \frac{1}{2 \sigma_{d_i}^2} \delta^2 .
    \end{align*}
    \begin{lemma}
        \label{lemma:lemma-2}
        Consider two vectors $\mathbf
        x$ and $\mathbf{y}$, such that $\mathbf{x} = \mathbf{y} + \mathbf{\Delta}$ and $\|\mathbf{\Delta}\| \ll \|\mathbf{y}\|$. Then $\frac{\mathbf{x}}{\|\mathbf{x}\|} - \frac{\mathbf{y}}{\|\mathbf{y}\|} \approx \frac{\mathbf{\Delta}}{\|\mathbf{y}\|}$.
    \end{lemma}
    \begin{proof}
        If $\|\mathbf{\Delta}\| \ll \|\mathbf{y}\|$, then $\|\mathbf{x}\| \approx \|\mathbf{y}\|$. So $\frac{\mathbf{x}}{\|\mathbf{x}\|} - \frac{\mathbf{y}}{\|\mathbf{y}\|} \approx \frac{\mathbf{y} + \mathbf{\Delta} - \mathbf{y}}{\|\mathbf{y}\|} = \frac{\mathbf{\Delta}}{\|\mathbf{y}\|}$.
    \end{proof}
    For the second term, we have
    \begin{align*}
        \mathbf{E}\left[ \log g_{u_i}(U_i;\theta, \rho) - \log g_{u_i}(U_i;\theta_0) \right] &= \mathbf{E}\left[ \kappa_i U_i^T \frac{\mathbf{z}}{\|\mathbf{z}\|} - \kappa_i U_i^T \frac{\mathbf{z}_0}{\|\mathbf{z}_0\|} \right] \\
        &\leq \kappa_i \left(\frac{\mathbf{z}_0}{\|\mathbf{z}_0\|}\right)^T \left( \frac{\mathbf{z}}{\|\mathbf{z}\|} - \frac{\mathbf{z}_0}{\|\mathbf{z}_0\|} \right) \\
        &\lesssim -\frac{\kappa_i}{\|\mathbf{z}_0\|} \delta .
    \end{align*}
    The first inequality comes from the fact that $\mathbf{E}(U_i) = \left( \coth(\kappa_i) - \frac{1}{\kappa_i} \right) \frac{\mathbf{z}_0}{\|\mathbf{z}_0\|} \preceq \frac{\mathbf{z}_0}{\|\mathbf{z}_0\|}$, for all $\kappa_i$~\cite{mardia2000directional}. The last inequality comes from the combination of Cauchy-Schwarz inequality and Lemmas~\ref{lemma:lemma-1} and~\ref{lemma:lemma-2}.
    \begin{lemma}
        \label{lemma:lemma-3}
        Consider the vectors $\mathbf{x}$ and $\mathbf{y}$, such that $\mathbf{x} = \mathbf{y} + \mathbf{\Delta}_1$ and $\|\mathbf{\Delta}\| \ll \|\mathbf{y}\|$, and the vectors $\mathbf{u}$ and $\mathbf{v}$, such that $\mathbf{u} = \mathbf{v} + \mathbf{\Delta}_2$. Then
        \begin{align*}
            \left( \frac{\mathbf{x}}{\|\mathbf{x}\|} \right)^T \mathbf{u} \ - \left( \frac{\mathbf{y}}{\|\mathbf{y}\|} \right)^T \mathbf{v} \ \approx \frac{\mathbf{y}^T \mathbf{\Delta}_2 + \mathbf{\Delta}_1^T \mathbf{\Delta}_2 + \mathbf{v}^T \mathbf{\Delta}_1}{\|\mathbf{y}\|}.
        \end{align*}
    \end{lemma}
    \begin{proof}
        The result is readily obtained from applying Lemma~\ref{lemma:lemma-2} and rearranging the expression.
    \end{proof}
    For the last term, let $M_i = \frac{2 f_{c,i}}{c}$, $w = \left(\frac{\mathbf{r} - \mathbf{s_i}}{\|\mathbf{r} - \mathbf{s_i}\|}\right)^T \mathbf{v}$ and $w_0 = \left(\frac{\mathbf{r}_0 - \mathbf{s_i}}{\|\mathbf{r}_0 - \mathbf{s_i}\|}\right)^T \mathbf{v}_0$, then we have
    \begin{align*}
        \mathbf{E}\left[ \log g_{f_i}(F_i;\theta, \rho) - \log g_{f_i}(F_i;\theta_0) \right] &= \mathbf{E}\left[ \frac{(F_i - M_i w_0)^2 - (F_i - M_i w)^2}{2 \sigma_{f_i}^2} \right] \\
        &= \mathbf{E}\left[ \frac{M_i^2 w_0^2 - 2 F_i M_i w_0 + 2 F_i M_i w - M_i^2 w^2}{2 \sigma_{f_i}^2} \right] \\
        &= -\frac{M_i^2 (w - w_0)^2}{2 \sigma_{f_i}^2} \\
        &\lesssim -\frac{\left( M_i \frac{\|\mathbf{z}_0\| + \|\mathbf{v}_0\| + 1}{\|\mathbf{z}_0\|} \right)^2}{2 \sigma_{f_i}^2} \delta^2 . \\
    \end{align*}
    The last inequality results from the combination Lemma~\ref{lemma:lemma-3} and Cauchy-Schwarz inequality. Combining the results for the three terms, we have
    \begin{align*}
        \sum_{i=1}^{n} \mathbf{E}\left[ \log g_i(X_i;\theta, \rho) - \log g_i(X_i;\theta_0) \right] \lesssim \sum_{i=1}^{n} \left( \frac{\delta^2}{2 \sigma_{d_i}^2} + \frac{\kappa_i}{\|\mathbf{z_0}\|} \delta + \frac{\left( M_i \frac{\|\mathbf{z}_0\| + \|\mathbf{v}_0\| + 1}{\|\mathbf{z}_0\|} \right)^2}{2 \sigma_{f_i}^2} \delta^2 \right) .
    \end{align*}
    Therefore, the inequality of the assumption is satisfied for
    \begin{align*}
        b_n = \frac{1}{2} \sum_{i=1}^{n} \left( \frac{\delta^2}{\sigma_{d_i}^2} + \frac{2 \kappa_i}{\|\mathbf{z_0}\|} \delta + \frac{\left( M_i \frac{\|\mathbf{z}_0\| + \|\mathbf{v}_0\| + 1}{\|\mathbf{z}_0\|} \right)^2}{\sigma_{f_i}^2} \delta^2 \right) .
    \end{align*}

    \item[(v)] For each $\theta \neq \theta_0$, and any sufficiently small values of $\rho$,
    \begin{align*}
        \sum_{i=1}^{\infty} \frac{1}{b_i^2} \mathbf{Var}\left[ \log g_i(X_i;\theta, \rho) - \log g_i(X_i;\theta_0) \right] < \infty .
    \end{align*}
    As before, choose $\delta_{\mathbf{r}}, \delta_{\mathbf{v}} > 0$ such that $\|\mathbf{r}' - \mathbf{r}\| \leq \delta_{\mathbf{r}}$, $\|\mathbf{v}' - \mathbf{v}\| \leq \delta_{\mathbf{v}}$ and $\|\mathbf{\theta}' - \mathbf{\theta}\| \leq \rho$, where $\mathbf{\theta} \coloneqq (\mathbf{r}, \mathbf{v})$. Let $\mathbf{z} = \mathbf{r} - \mathbf{s}_i$ and $\mathbf{z}_0 = \mathbf{r}_0 - \mathbf{s}_i$. Given the constraints on $\delta$ and $\rho$, by the triangle inequality let $\mathbf{r} = \mathbf{r}_0 + \mathbf{\Delta_{r}}$, such that $\|\mathbf{\Delta_{r}}\| \leq \delta_{\mathbf{r}}$, and $\mathbf{v} = \mathbf{v}_0 + \mathbf{\Delta_{v}}$, such that $\|\mathbf{\Delta_{v}}\| \leq \delta_{\mathbf{v}}$. Let $\delta = \max\{\delta_{\mathbf{r}}, \delta_{\mathbf{v}}\}$. Due to the assumption of independence between measurements, note that
    \begin{align*}
        \mathbf{Var}\left[ \log g_i(X_i;\theta, \rho) - \log g_i(X_i;\theta_0) \right] &= \mathbf{Var}\left[ \log g_{d_i}(D_i;\theta, \rho) - \log g_{d_i}(D_i;\theta_0) \right] \\
        &+ \mathbf{Var}\left[ \log g_{u_i}(U_i;\theta, \rho) - \log g_{u_i}(U_i;\theta_0) \right] \\
        &+ \mathbf{Var}\left[ \log g_{f_i}(F_i;\theta, \rho) - \log g_{f_i}(F_i;\theta_0) \right] .
    \end{align*}
    So for the first term,
    \begin{align*}
        \mathbf{Var}\left[ \log g_{d_i}(D_i;\theta, \rho) - \log g_{d_i}(D_i;\theta_0) \right] &= \mathbf{Var}\left[ \frac{(D_i - \|\mathbf{z}_0\|)^2 - (D_i - \|\mathbf{z}\|)^2}{2 \sigma_{d_i}^2} \right] \\
        &= \mathbf{Var}\left[ \frac{\|\mathbf{z}_0\|^2 - 2 D_i \|\mathbf{z}_0\| + 2 D_i \|\mathbf{z}\| - \|\mathbf{z}\|^2}{2 \sigma_{d_i}^2} \right] \\
        &= \frac{(\|\mathbf{z}\| - \|\mathbf{z}_0\|)^2}{\sigma_{d_i}^4} \mathbf{Var}\left[ D_i \right] = \frac{\delta^2}{\sigma_{d_i}^2} .
    \end{align*}
    For the second term, we have
    \begin{align*}
        \mathbf{Var}\left[ \log g_{u_i}(U_i;\theta, \rho) - \log g_{u_i}(U_i;\theta_0) \right] &= \mathbf{Var} \left[ \kappa_i U_i^T \frac{\mathbf{z}}{\|\mathbf{z}\|} - \kappa_i U_i^T \frac{\mathbf{z}_0}{\|\mathbf{z}_0\|} \right] \\
        &\approx \kappa_i^2 \mathbf{Var} \left[ U_i^T \frac{\mathbf{\Delta_r}}{\|\mathbf{z}_0\|} \right] \\
        &= \frac{\kappa_i^2}{\|\mathbf{z}_0\|^2} \mathbf{\Delta_r}^T \mathbf{Var} \left[ U_i \right] \mathbf{\Delta_r} \\
        &\leq \frac{\kappa_i^2}{\|\mathbf{z}_0\|^2} \mathbf{\Delta_r}^T \left(\mathbf{I} + \frac{\mathbf{z}_0}{\|\mathbf{z}_0\|} \left(\frac{\mathbf{z}_0}{\|\mathbf{z}_0\|}\right)^T \right) \mathbf{\Delta_r} \\
        &= \frac{2 \delta^2 \kappa_i^2}{\|\mathbf{z}_0\|^2} \\
        &\leq \frac{2 \delta \kappa_i}{\|\mathbf{z}_0\|}, \quad\quad \mbox{for sufficiently small } \delta .
    \end{align*}
    The first inequality results from the fact that $\mathbf{Var}(U_i) \preceq \mathbf{I} + \frac{\mathbf{z}_0}{\|\mathbf{z}_0\|} \left(\frac{\mathbf{z}_0}{\|\mathbf{z}_0\|} \right)^T$. The last inequality holds for sufficiently small $\delta$. For the last term, we obtain
    \begin{align*}
        \mathbf{Var}\left[ \log g_{f_i}(F_i;\theta, \rho) - \log g_{f_i}(F_i;\theta_0) \right] &= \mathbf{Var}\left[ \frac{(F_i - M_i w_0)^2 - (F_i - M_i w)^2}{2 \sigma_{f_i}^2} \right] \\
        &= \mathbf{Var}\left[ \frac{M_i^2 w_0^2 - 2 F_i M_i w_0 + 2 F_i M_i w - M_i^2 w^2}{2 \sigma_{f_i}^2} \right] \\
        &= \frac{M_i^2 (w - w_0)^2}{ \sigma_{f_i}^4} \mathbf{Var}\left[ F_i \right]  \\
        &\lesssim \frac{\left( M_i \frac{\|\mathbf{z}_0\| + \|\mathbf{v}_0\| + 1}{\|\mathbf{z}_0\|} \right)^2}{\sigma_{f_i}^2} \delta^2 . \\
    \end{align*}
    Combining the results for the three terms, we have
    \begin{align*}
        \mathbf{Var}\left[ \log g_i(X_i;\theta, \rho) - \log g_i(X_i;\theta_0) \right] \lesssim \left( \frac{\delta^2}{\sigma_{d_i}^2} + \frac{2 \kappa_i}{\|\mathbf{z_0}\|} \delta + \frac{\left( M_i \frac{\|\mathbf{z}_0\| + \|\mathbf{v}_0\| + 1}{\|\mathbf{z}_0\|} \right)^2}{\sigma_{f_i}^2} \delta^2 \right) .
    \end{align*}
    So, if we remember $b_n$ in the last assumption, with this result, we see that 
    \begin{align*}
        b_i = \frac{1}{2} \sum_{j=1}^{i} \mathbf{Var}\left[ \log g_j(X_i;\theta, \rho) - \log g_j(X_i;\theta_0) \right] .
    \end{align*} 
    Therefore, the inequality of this assumption holds due to the following lemma:
    \begin{lemma}
        Let $a_n$ be a sequence of non-negative real numbers, and let $S_n = a_1 + a_2 + \dots + a_n$. Then $\sum_{n=1}^{\infty} \frac{a_n}{S_n^2} < \infty$.
    \end{lemma}
\end{itemize}

So, we finalize our proof and show that the Maximum Likelihood Estimator for orbit determination in~\cite{ferreira2024generalizing} is strongly consistent, i.e., the estimated parameters converge almost surely to the true parameters, as it fulfills the necessary assumptions.

%% file: main.bbl
\begin{thebibliography}{10}

\bibitem{wald1949note}
A.~Wald, ``Note on the consistency of the maximum likelihood estimate,'' {\em The Annals of Mathematical Statistics}, vol.~20, no.~4, pp.~595--601, 1949.

\bibitem{goel1974note}
P.~K. Goel, ``A note on the consistency of maximum likelihood estimators,'' {\em Scandinavian Actuarial Journal}, vol.~1974, no.~4, pp.~211--220, 1974.

\bibitem{ferguson2017course}
T.~S. Ferguson, {\em A Course in Large Sample Theory}.
\newblock Routledge, 1~ed., 2017.

\bibitem{hald1999history}
A.~Hald, ``On the history of maximum likelihood in relation to inverse probability and least squares,'' {\em Statistical Science}, vol.~14, no.~2, pp.~214--222, 1999.

\bibitem{bishop2006pattern}
C.~M. Bishop, {\em Pattern recognition and machine learning}, vol.~4.
\newblock Springer.

\bibitem{rao1973linear}
C.~R. Rao, C.~R. Rao, M.~Statistiker, C.~R. Rao, and C.~R. Rao, {\em Linear statistical inference and its applications}, vol.~2.
\newblock Wiley New York, 1973.

\bibitem{ferreira2024generalizing}
R.~Ferreira, F.~Valdeira, M.~Guimar{\~a}es, and C.~Soares, ``Generalizing trilateration: Approximate maximum likelihood estimator for initial orbit determination in low-earth orbit,'' 2024.

\bibitem{sen1994large}
P.~K. Sen and J.~M. Singer, {\em Large Sample Methods in Statistics: An Introduction with Applications}, vol.~25.
\newblock CRC press, 1994.

\bibitem{rudin1964principles}
W.~Rudin {\em et~al.}, {\em Principles of mathematical analysis}, vol.~3.
\newblock McGraw-hill New York, 1964.

\bibitem{mardia2000directional}
K.~V. Mardia, P.~E. Jupp, and K.~Mardia, {\em Directional statistics}, vol.~2.
\newblock Wiley Online Library, 2000.

\end{thebibliography}
